\documentclass{ifacconf}
\usepackage{graphicx}      
\usepackage{enumitem}
\usepackage{siunitx}
\usepackage{tikz}
\usepackage{comment}
\usepackage[tbtags]{mathtools}
\usepackage{diagbox}
\usepackage{mathrsfs}
\usepackage{amssymb,amsfonts}
\usepackage{algorithmic}
\mathtoolsset{showonlyrefs=true}
\usepackage[dvipsnames]{xcolor}

\usepackage{float}
\usepackage{accents}
\usepackage{multicol}
\usepackage{array}[=2023/01/01]
\usepackage{multirow}
\usepackage{hhline}
\usepackage{arydshln}
\usepackage{bbm}
\usepackage{bm}
\usepackage{enumitem}
\usepackage{subcaption}

\usepackage{soul}
\DeclareMathOperator{\im}{im}

\usepackage{natbib}

\newcommand{\R}{\mathbb{R}}

\begin{document}
\begin{frontmatter}

\title{Data-driven Synchronization for Network Systems with Noiseless Data\thanksref{footnoteinfo}} 

\thanks[footnoteinfo]{The work of Yongzhang Li was supported by the China Scholarship Council under Grant 202106070023.}

\author[First]{Yongzhang Li} 
\author[First]{M. K. Camlibel}

\address[First]{Bernoulli Institute for Mathematics, Computer Science and Artificial Intelligence, University of Groningen, The Netherlands, (e-mail: yongzhang.li@rug.nl, m.k.camlibel@rug.nl).}

\begin{abstract}                
For a collection of homogeneous LTI systems that is interconnected by a protocol, given the network topology and the system model, one may obtain a feedback gain to synchronize the network. However, the model-based methods cannot be applied in case the system model is \emph{unknown}. Therefore, in this paper, we study the data-driven synchronization problem for homogeneous networks. In particular, given a collection of LTI systems, we collect the input-state data from one individual system. Then, given the network topology, we provide data-based necessary and sufficient conditions for synchronizability. Once the conditions are satisfied, one can also obtain a feedback gain directly from data to synchronize the network with the corresponding design method provided in this paper. Finally, we illustrate our results with a numerical simulation.
\end{abstract}

\end{frontmatter}

\section{Introduction}
Over the past two decades, the synchronization problem for network systems has been an active topic in the research of systems and control. For a collection of individual systems, the design of an interconnection protocol to achieve synchronization is an important topic in many fields, e.g., multi-agent systems, social networks, and biological networks (see \cite{sivrikaya2004time,arenas2008synchronization,proskurnikov2015opinion,dong2018survey,saberi2022syn}). The traditional methods of studying synchronization problems are model-based. In particular, employing the conditions for synchronizability and corresponding design methods require the individual system models to be \emph{priori known} and accurate enough. However, precise system models are not always available.

In order to avoid the dependency on system models, some data-based methods were developed to study synchronization problems. For example, reinforcement learning and iteration algorithms have gained popularity over the last decade. In \cite{Li2017}, and later in \cite{Qin2019}, methods based on off-policy reinforcement learning were proposed to realize optimal synchronization for homogeneous multi-agent systems. Additionally, one may be interested in a two-stage policy iteration algorithm proposed in \cite{peng2019data} for solving the optimal tracking problem in leader-follower multi-agent systems, which builds upon the traditional policy iteration algorithms of adaptive dynamic programming. These methods enable one-step design from data to bypass the lack of knowledge about system models, but typically require a huge amount of real-time data. As such, these learning-based processes, especially for network systems, can meet challenges in data collection, transportation, and processing. This makes it interesting to study data-based methods that are less computationally expensive.

In recent years, data-driven control has received considerable attention. We refer the reader to the recent book \cite{DBLSCT} for a comprehensive overview that includes the historical background (see Chapter 1) and the current state of the art. In this paper, we adopt the \emph{data informativity} framework introduced in \cite{van2020data}. This framework allows system analysis and control directly from offline measured data, without the requirement of persistency of excitation. Several fundamental analysis and control problems were already studied in \cite{van2020data}. Later, this framework was extended to a noisy data setting, and several control problems were addressed based on matrix versions of S-Lemma in \cite{van2020noisy,van2022data,HenkQMI2023,van2023inf}.


With the data informativity framework and the development of corresponding tools, several studies have been done on synchronization problems. A study on output synchronization for heterogeneous leader-follower multi-agent systems has been conducted in \cite{jiao2021data}, which provided a direct design method for feedback gains for each follower. Also, a data-based distributed control method was proposed in \cite{eising2022informativity} to synchronize heterogeneous network systems with a given sparsity structure. Furthermore, they minimized the sparsity structure for the network. These methods, distinct from learning-based or iteration-based methods, rely on only offline measured data. See also \cite{steentjes2021h,celi2023distributed} for other related works. The proposed works presented reliable methods for designing the network structure and controllers to synchronize the heterogeneous network systems.

Different from the works on heterogeneous network systems mentioned above, in this paper, we study the data-driven synchronization for \emph{homogeneous} network systems. In homogeneous network systems, all individual systems are identical. Therefore, in a noiseless data-driven setting, we collect data from only one system, but analyze and control the whole network. For a given network topology, we provide necessary and sufficient conditions on measured data such that a feedback gain can be obtained for an interconnection protocol to achieve synchronization with the corresponding design method.

The remainder of the paper is organized as follows. In Section \ref{Ch1_Sec_Syn_Mod}, we introduce the synchronization problem for homogeneous network systems and provide model-based results. In Section \ref{Ch1_Sec_Syn_dd}, we formulate and study the data-driven synchronization problems with offline measured data. A numerical simulation is provided in Section \ref{Ch1_Sec_Syn_ex} to illustrate the results. Finally, Section \ref{Ch1_Sec_conclusion} concludes the paper.  

\subsection*{Notation}
The $p$-dimensional vector whose elements are all $1$ is denoted by $\mathbbm{1}_p$. Given a positive integer $p$, we define $\bar{p}=\{1,2,\dots , p\}$. The set of $n\times n$ real symmetric matrices is denoted by $\mathbb{S}^n$. We say matrix $A$ is positive definite (resp., positive semi-definite), denoted by $A>0$ (resp., $A\geq 0$) if $A\in\mathbb{S}^{n}$ and all its eigenvalues are positive (resp., nonnegative). We say matrix $A$ is negative definite (resp., negative semi-definite), denoted by $A<0$ (resp., $A\leq 0$) if $-A>0$ (resp., \mbox{$-A\geq 0$}). A right inverse of a real matrix $A$ is denoted by $A^R$. The Kronecker product of real matrices $A$ and $B$ is denoted by $A \otimes B$. The set of all eigenvalues of a real square matrix $A$ is denoted by $\sigma(A)$. A square matrix is called \emph{Schur} if all its eigenvalues lie inside the open unit disc in the complex plane.

\section{Synchronization for network systems}
\label{Ch1_Sec_Syn_Mod}
We consider a collection of $p$ discrete-time LTI systems of the form
\begin{equation}
\label{eq:Ch1_mod_LTI}
    \bm{x}_i(k+1)=A\bm{x}_i(k)+B\bm{u}_i(k),
\end{equation}
where $\bm{x}_i\in \R^n$ and $\bm{u}_i\in\R^m$ are, respectively, the state and the input of the $i$th system with $i \in \bar{p}$. The matrices $A$ and $B$ are of suitable sizes. The systems \eqref{eq:Ch1_mod_LTI} are interconnected by a protocol of the form
\begin{equation}
\label{eq:Ch1_intcont_protocol}
    \bm{u}_i(k)=K\sum_{j = 1}^{p} {c_{ij} \bm{x}_j(k)}
\end{equation}
where $K \in \R^{m \times n}$ and $c_{ij} \in \R$ for $i,j \in \bar{p}$. Now, let
\begin{equation}
\bm{x}=\begin{bmatrix} \bm{x}_1 \\ \bm{x}_2 \\ \vdots \\ \bm{x}_p \end{bmatrix}, \
C=\begin{bmatrix} c_{11} & c_{12} & \cdots & c_{1p} \\ c_{21} & c_{22} & \cdots &
	c_{2p} \\ \vdots & \vdots & \ddots & \vdots \\ c_{p1} & c_{p2} & \cdots & c_{pp}\end{bmatrix}.
\end{equation} Then, the overall network \eqref{eq:Ch1_mod_LTI}-\eqref{eq:Ch1_intcont_protocol} can be more compactly described by
\begin{equation}
\bm{x}(k+1)=(I_p\otimes A+C\otimes BK)\bm{x}(k).
\label{eq:Ch1_overall_dy}
\end{equation}
We call the matrix $C$ the \emph{interconnection matrix}, which represents the network topology of the network \eqref{eq:Ch1_overall_dy}.

Next, we define the synchronization for network systems as follows.
\begin{defn}
    The network \eqref{eq:Ch1_overall_dy} is said to be synchronized if for every trajectory $\bm{x}$ of \eqref{eq:Ch1_overall_dy} and for all $i,j \in \bar{p}$, $\lim\limits_{k \rightarrow \infty} \bm{x}_i(k)-\bm{x}_j(k) = 0$.
\end{defn}
We have the following blanket assumption throughout the paper. 
\begin{assum}
\label{ass:Ch1_assC_undirect}
    The vector $\mathbbm{1}_p$ is the eigenvector of $C$ corresponding to the real eigenvalue $\mu$ whose algebraic multiplicity is $1$, and all other eigenvalues of $C$ are real and positive. 
\end{assum}
There are various common choices for the interconnection matrix $C$ that satisfy Assumption \ref{ass:Ch1_assC_undirect}. For example, in multi-agent systems with an undirected communication network, the Laplacian matrix of the graph that represents the network topology is a common choice (see \cite{you2011}). One can also find that symmetric row stochastic matrices are used to represent the network systems (see \cite{zhang2015syn}). Furthermore, the left random-walk Laplacian matrix of a simple graph is also commonly used if a normalized protocol is required (see \cite{saberi2022syn}).


The following proposition provides necessary and sufficient conditions for synchronization.
\begin{prop}
   The network system \eqref{eq:Ch1_overall_dy} is synchronized if and only if $A+\lambda BK$ is Schur for all $\lambda \in \sigma(C)\backslash\{\mu\}$. 
\end{prop}

\begin{pf}
    Let $\mathcal{L}$ be the Laplacian matrix of a complete undirected graph with $p$ vertices. We define 
    \begin{equation}
    \label{eq:Ch1_output_L}
        \bm{y}(k)\coloneqq(\mathcal{L} \otimes I_n)\bm{x}(k).
    \end{equation}
    One can verify that $\lim\limits_{k \rightarrow \infty}\bm{y}(k)= 0$ for every trajectory $\bm{x}$ if and only if $\lim\limits_{k \rightarrow \infty} \bm{x}_i(k)-\bm{x}_j(k) = 0$ for all $i,j \in \bar{p}$. Next, we prove that $\lim\limits_{k \rightarrow \infty}\bm{y}(k)= 0$ for every trajectory $\bm{x}$ if and only if $A+\lambda BK$ is Schur for all $\lambda \in \sigma(C)\backslash\{\mu\}$. To this end, let $\lambda_i$ be the $i$th eigenvalue of $C$ except $\mu$. Since $C$ is a real square matrix with only real eigenvalues, it follows from Schur's theorem (see e.g.\cite[Sec. 2.3]{horn2012matrix}) that there exists an orthogonal matrix $V \in \R^{p \times p}$ such that 
    \begin{equation}
        \Lambda_C\coloneqq V^\top CV=\begin{bmatrix}
                 \mu & *          & \cdots & *      \\
                 0   & \lambda_1  & \cdots & *      \\
                 0   & 0          & \ddots & \vdots \\
                 0   & 0          & \cdots & \lambda_{p-1}
                 \end{bmatrix}.
    \end{equation} 
    Now, define
    \begin{equation}
        \bm{z}\coloneqq (V\otimes I_n)\bm{x}.
    \end{equation}
    From \eqref{eq:Ch1_overall_dy} and \eqref{eq:Ch1_output_L}, we have
    \begin{equation}
        \bm{z}(k+1)=((I_p \otimes A)+(\Lambda_C\otimes BK))\bm{z}(k),
    \end{equation}
    and
    \begin{equation}
        \bm{y}(k+1)=(\mathcal{L}V^\top\otimes I_n)((I_p \otimes A)+(\Lambda_C\otimes BK))\bm{z}(k).
    \end{equation}
    According to Assumption \ref{ass:Ch1_assC_undirect}, we have $\mathbbm{1}_p$ as an eigenvector of $C$ corresponding to the eigenvalue $\mu$. This implies that
    \begin{equation}
        V=\begin{bmatrix}
            \frac{\mathbbm{1}_p}{\sqrt{p}} & V_1
        \end{bmatrix}.
    \end{equation}
    Now, let $A_K(\lambda)\coloneqq A+\lambda BK$. Since $\mathbbm{1}_p$ is also an eigenvector of $\mathcal{L}$ corresponding to the eigenvalue $0$, for any initial state $\bm{z}(0)$, we have 
    \begin{equation}
    \begin{split}
        &\bm{y}(k+1)=\\
        &\left(\begin{bmatrix}
            0 &  \mathcal{L}V_1
        \end{bmatrix}\otimes I_n\right)\begin{bmatrix}
            A_K^k(\mu) & *              & \cdots & * \\
            0        & A_K^k(\lambda_1) & \cdots & * \\
            0        & 0              & \ddots & \vdots \\
            0        & 0              & \cdots & A_K^k(\lambda_{p-1})
        \end{bmatrix}\bm{z}(0). 
    \end{split}
    \end{equation}
    Since $\mathcal{L}V_1$ has full column rank, $\lim\limits_{k \rightarrow \infty}\bm{y}(k)= 0$ for every trajectory $\bm{z}(0)$ if and only if $A+\lambda BK$ is Schur for all $\lambda \in \sigma(C)/\{\mu\}$. This completes the proof. \hfill \qed
\end{pf}

This proposition provides conditions for a \textit{given} network to be synchronized. Another interesting question is how to design a gain $K$ for given $(A,B)$ and $C$ to synchronize the network system. An answer to this question is given by the following proposition, where $\lambda_{\text{min}}$  ($\lambda_{\text{max}}$) denotes the smallest (largest) element of the set $\sigma(C)\backslash\{\mu\}$.

	\begin{prop}
		Let $Q>0$ and $R>0$. Suppose that $P>0$ satisfies the algebraic Riccati equation
		\begin{equation}
		P-A^\top PA+\frac{\lambda_{\text{min}}}{\lambda_\text{max}} A^\top PB(B^\top PB+R)^{-1}B^\top PA-Q=0.
		\label{Riccati_pos}
		\end{equation} Then, the network \eqref{eq:Ch1_overall_dy} is synchronized by the choice 
		\begin{equation}
		K=-\frac{1}{\lambda_\text{max}}(B^\top PB+R)^{-1}B^\top PA.
		\label{design_gain}
		\end{equation}
		\label{thm:Ch1_riccati}
	\end{prop}

\begin{pf}
   Let $H_\lambda:=P-(A+\lambda BK)^\top P (A+\lambda BK)$. We observe that  $A+\lambda BK$ is Schur whenever $H_\lambda > 0$. Therefore, it is enough to show that $H_\lambda > 0$ for all $\lambda \in [\lambda_{\text{min}},\lambda_{\text{max}}]$ and hence for all $\lambda \in \sigma(C)\backslash\{\mu\}$.

   Let $\lambda \in [\lambda_{\text{min}},\lambda_{\text{max}}]$. By writing the expansion of $H_\lambda$, we have
   \begin{align}
        H_\lambda &= P - A^\top PA-\lambda A^\top PBK-\lambda K^\top B^\top PA   \\
          & \ \ \ -\lambda^2K^\top B^\top PBK \\
          &\overset{\text{\eqref{design_gain}}}{=} P-A^\top PA+2\lambda \lambda_{\text{max}} K^\top (B^\top PB+R) K \\
          & \ \ \ -\lambda^2K^\top B^\top PBK \\
          &\overset{R>0}{>} P-A^\top PA+(2\lambda \lambda_{\text{max}}-\lambda^2) K^\top (B^\top PB+R)K,
   \end{align} where the second equality is obtained by replacing $A^\top PB$ and $B^\top PA$  with $-\lambda_{\text{max}}K^\top (B^\top PB+R)$ and $-\lambda_{\text{max}}(B^\top PB+R)K$ respectively based on \eqref{design_gain}.
   
   Then, it follows from \eqref{Riccati_pos} that 
	\begin{equation}
	\begin{split}
	H_\lambda&>Q+(2\lambda \lambda_{\text{max}}-\lambda^2-\lambda_{\text{min}} \lambda_{\text{max}}) K^\top (B^\top PB+R)K \\
	 &>Q+(\lambda \lambda_{\text{max}}-\lambda^2)K^\top (B^\top PB+R)K.
	\end{split}
    \label{ineq_H}
	\end{equation}
	
	Since $Q>0$, $R>0$, and $P>0$, the inequality \eqref{ineq_H} implies that $H_\lambda>0$ for all $\lambda \in [\lambda_{\text{min}},\lambda_{\text{max}}]$. This proves the proposition. \hfill \qed
\end{pf}

\section{Data-Driven synchronization}
\label{Ch1_Sec_Syn_dd}
In Section \ref{Ch1_Sec_Syn_Mod}, we introduced the synchronization problem for network systems and provided necessary and sufficient conditions for synchronization in a model-based setting. In this section, we turn our attention to the data-driven setting.

Consider a discrete-time LTI system of the form
\begin{equation}
\label{eq:Ch1_mod_true}
    x(k+1)=A_\textup{true}x(k)+B_\textup{true}u(k),
\end{equation}
where $x(k) \in \R^n$ and $u \in \R^m$ are, respectively, the state and input. We refer to this system as the \emph{true system}.

We assume that the matrices $A_\textup{true}$ and $B_\textup{true}$ are unknown but we collect the input-state data
\begin{equation}
    \mathcal{D}\coloneqq \left(\begin{bmatrix}
    u(0)  & \cdots & u(T-1)
\end{bmatrix}, \begin{bmatrix}
    x(0) & \cdots & x(k)
\end{bmatrix}\right)
\end{equation}
from the system \eqref{eq:Ch1_mod_true}. We also define the data matrices
\begin{equation}
\begin{split}
U_-&\coloneqq\begin{bmatrix} u(0) & u(1) & \cdots & u(T-1) \end{bmatrix},  \\
X_-&\coloneqq\begin{bmatrix} x(0) & x(1) & \cdots & x(T-1) \end{bmatrix}, \\
X_+&\coloneqq\begin{bmatrix} x(1) & x(2) & \cdots & x(k) \end{bmatrix}.
\end{split}
\end{equation}
A pair of real matrices $(A,B)$ is called a \emph{data-consistent system} if it satisfies
\begin{equation}
    X_+=AX_-+BU_-.
\end{equation} 
We define the set of all data-consistent systems as 
\begin{equation}
\label{eq:Ch1_Sigma_D}
    \Sigma_\mathcal{D}\coloneqq\{(A,B):X_+=AX_-+BU_-\}.
\end{equation} 
Since the data $\mathcal{D}$ are collected from \eqref{eq:Ch1_mod_true}, clearly we have $(A_\textup{true},B_\textup{true}) \in \Sigma_\mathcal{D}$. Then, to ensure the network \eqref{eq:Ch1_overall_dy} with $(A,B)=(A_\textup{true},B_\textup{true})$ is synchronized, we have to guarantee that \eqref{eq:Ch1_overall_dy} is synchronized for all \mbox{$(A,B) \in \Sigma_\mathcal{D}$}. Following the idea of the informativity framework proposed in \cite{van2020data,van2023inf,DBLSCT}, we provide the following definition for the data-based synchronizability notion.
\begin{defn}
    The data $\mathcal{D}$ are said to be \emph{informative for synchronizability} with respect to the interconnection matrix $C$ if there exists $K$ such that $A+\lambda BK$ is Schur for all $\lambda \in \sigma(C)/\{\mu\}$ and for all $(A,B) \in \Sigma_\mathcal{D}$.
\end{defn}

Let $\sigma(C)\backslash\{\mu\}=\{\lambda_1,\lambda_2,\dots,\lambda_\ell\}$. Next, we provide the necessary and sufficient conditions for the data $\mathcal{D}$ to be informative for synchronizability with respect to a given interconnection matrix $C$.
\begin{thm}
\label{thm:Ch1_DDsyn_noiseless}
    The data $\mathcal{D}$ are informative for synchronizability with respect to the interconnection matrix $C$ if and only if the following conditions are satisfied:
    \begin{enumerate}[label=(\alph*)]
        \item $X_-$ has full row rank.
        \item There exist right inverses $(X_-^R)_1,\dots,(X_-^R)_{\ell}$ of $X_-$ s.t. $X_+(X_-^R)_i$ is Schur and
              \begin{equation}
              \label{eq:Ch1_data_Rinverse}
              U_-(\lambda_i(X_-^R)_j-\lambda_j(X_-^R)_i)=0
              \end{equation}
              for all $i,j \in \bar \ell$.
    \end{enumerate}
    Moreover, if such right inverses of $X_-$ exist, then the network \eqref{eq:Ch1_overall_dy} with $(A,B)\in\Sigma_\mathcal{D}$ is synchronized by the choice $K=U_-(X_-^R)_1/\lambda_1$. 
\end{thm}

To prove this theorem, we need the following auxiliary result.

\begin{lem}
\label{lem:Ch1_DDsyn_im_K}
    Suppose that the data $\mathcal{D}$ are informative for synchronizability with respect to the interconnection matrix $C$, and $K$ is a gain such that $A+\lambda BK$ is Schur for all $\lambda \in \sigma(C)/\{\mu\}$ and all $(A,B) \in \Sigma_\mathcal{D}$. Then, we have 
    \begin{equation}
        \im \begin{bmatrix} I\\\lambda K \end{bmatrix} \subseteq \im \begin{bmatrix} X_- \\ U_- \end{bmatrix}
    \end{equation}
    for all $\lambda \in \sigma(C)/\{\mu\}$.
\end{lem}

\begin{pf}
    We define
    \begin{equation}
    \label{eq:Ch1_Sigma_0}
        \Sigma_0\coloneqq\{(A_0,B_0) :A_0X_-+B_0U_-=0\}.
    \end{equation}
    Let $\lambda \in \sigma(C)/\{\mu\} $. We first claim that $A_0+\lambda B_0K$ is \emph{nilpotent} for all $(A_0,B_0) \in \Sigma_0$. To prove this claim, let $A_\alpha= A+\alpha A_0$ and $B_\alpha= B+\alpha B_0$ for some $\alpha \in \R$. Then, it follows from \eqref{eq:Ch1_Sigma_D} and \eqref{eq:Ch1_Sigma_0} that $(A_\alpha,B_\alpha) \in \Sigma_{\mathcal{D}}$ for all $\alpha$, all $(A,B) \in \Sigma_\mathcal{D}$, and all $(A_0,B_0) \in \Sigma_0$. Since $A+\lambda BK$ is Schur for all $(A,B) \in \Sigma_\mathcal{D}$, we observe that 
    \begin{equation}
        A+\lambda BK+\alpha(A_0+\lambda B_0K)
    \end{equation}
    is Schur for all $\alpha$, all $(A,B) \in \Sigma_\mathcal{D}$ and all $(A_0,B_0) \in \Sigma_0$. This implies that $A_0+\lambda B_0K$ is nilpotent for all \mbox{$(A_0,B_0) \in \Sigma_0$}. This proves the claim. Next, note that we have
    \begin{equation}
        ((A_0+\lambda B_0K)^\top A_0,(A_0+\lambda B_0K)^\top B_0) \in \Sigma_0
    \end{equation}
    whenever $(A_0,B_0) \in \Sigma_0$. This implies that $$(A_0+\lambda B_0K)^\top(A_0+\lambda B_0K)$$ is nilpotent. Since the only symmetric nilpotent matrix is the zero matrix, we conclude that $A_0+\lambda B_0 K=0$ for all $(A_0,B_0) \in \Sigma_0$. This is equivalent to
    \begin{equation}
        \ker \begin{bmatrix} X_-^\top & U_-^\top \end{bmatrix} \subseteq \ker \begin{bmatrix} I&\lambda K^\top \end{bmatrix},
    \end{equation}
    which is equivalent to
    \begin{equation}
        \im \begin{bmatrix} I\\\lambda K \end{bmatrix} \subseteq \im \begin{bmatrix} X_- \\ U_- \end{bmatrix}.
    \end{equation}
    This completes the proof. \hfill \qed
\end{pf}

Now, we are in a position to prove Theorem \ref{thm:Ch1_DDsyn_noiseless}.

\textbf{Proof of Theorem~\ref{thm:Ch1_DDsyn_noiseless}:}
We first prove the ``if' part. We observe from \eqref{eq:Ch1_data_Rinverse} that
\begin{equation}
    \frac{U_-(X_-^R)_i}{\lambda_i}=\frac{U_-(X_-^R)_j}{\lambda_j}
\end{equation}
for all $i,j \in \bar \ell$. This implies that  $U_-(X_-^R)_i/\lambda_i$ is the same for all $i \in \bar \ell$. We take $K_i=U_-(X_-^R)_i/\lambda_i$, and have 
\begin{equation}
\label{eq:Ch1_ABK_eq_X+X-}
     A+\lambda_i BK_i=(AX_-+BU_-)(X_-^R)_i=X_+(X_-^R)_i.
\end{equation}
This, together with the fact that $X_+(X_-^R)_i$ is Schur for all $i \in \bar \ell$, implies that $A+\lambda BK$ is Schur with $K=K_1$ for all $\lambda \in \sigma(C)/\{\mu\}$ and all $(A,B) \in \Sigma_{\mathcal{D}}$.

Next, we prove the ``only if'' part. Let $K$ be such that $A+\lambda BK$ is Schur with $K$ for all $\lambda \in \sigma(C)/\{\mu\}$ and all $(A,B) \in \Sigma_\mathcal{D}$. According to Lemma \ref{lem:Ch1_DDsyn_im_K}, we have
\begin{equation}
    \im \begin{bmatrix} I\\\lambda K \end{bmatrix} \subseteq \im \begin{bmatrix} X_- \\ U_- \end{bmatrix}
\end{equation}
for all $\lambda \in \sigma(C)/\{\mu\}$. This implies that $X_-$ has full row rank, and there exists a sequence of right inverses $(X_-^R)_1,\dots,(X_-^R)_{\ell}$ of $X_-$ such that 
\begin{equation}
\label{eq:Ch1_IK_equal_data}
    \begin{bmatrix} I\\\lambda_iK \end{bmatrix}=\begin{bmatrix} X_- \\ U_- \end{bmatrix}(X_-^R)_i
\end{equation}
for all $i \in \bar \ell$. Now, we observe from \eqref{eq:Ch1_IK_equal_data} that $$K=U_-(X_-^R)_i/\lambda_i=U_-(X_-^R)_j/\lambda_j$$ for all $i,j \in \bar \ell$. This implies that 
\begin{equation}
    U_-(\lambda_i(X_-^R)_j-\lambda_j(X_-^R)_i)=0
\end{equation}
for all $i,j \in \bar \ell$. Note that we have
\begin{equation}
\label{eq:Ch1_ABK_eq_X+X-_reverse}
     X_+(X_-^R)_i=(AX_-+BU_-)(X_-^R)_i=A+\lambda_i BK.
\end{equation}
This implies that $X_+(X_-^R)_i$ is Schur for all $i \in \bar \ell$. This completes the proof.\hfill \qed

\section{Illustrative Example}
\label{Ch1_Sec_Syn_ex}
\begin{figure}[h!]
    \centering
    \includegraphics[width=0.5\linewidth]{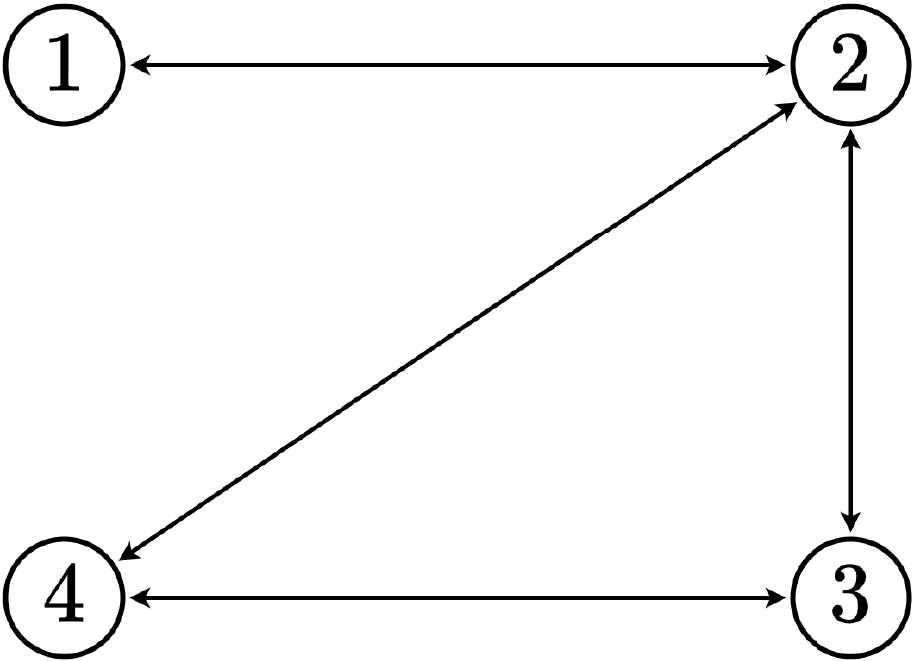}
    \caption{The network topology represented by $C$ as in \eqref{eq:Ch1_ex_noiseless_C}.}
    \label{fig:Ch1_ex_noiseless_topo}
\end{figure}

In this section, we provide a numerical simulation to illustrate the result we provided in Theorem \ref{thm:Ch1_DDsyn_noiseless}. Consider the true system given by
\begin{equation}
\label{eq:Ch1_ex_truesystem}
    A_\text{true}=\begin{bmatrix}
       \ 1 \ & 1 \ & 0 \ \\ \ 0 \ & 1 \ & 1 \ \\ \ 0 \ & 0 \ & 1 \
    \end{bmatrix}, B_\text{true}=\begin{bmatrix}
       \ 1 \ & 1 \ \\ \ 1 \ & 0 \ \\ \ 0 \ & 1 \
    \end{bmatrix}.
\end{equation}
Suppose that $p=4$ and the interconnection matrix is given by
\begin{equation}
\label{eq:Ch1_ex_noiseless_C}
    C=\begin{bmatrix}
      \  1 \ & -1 \ & 0 \ & 0 \ \\
      \  -1 \ & 3 \ & -1 \ & -1 \ \\
      \ 0 \ & -1 \ & 2 \ & -1 \ \\
      \ 0 \ & -1 \ & -1 \ & 2 \ \\
    \end{bmatrix}.
\end{equation}
For this $C$, we have $\mu=0$ and $\sigma(C)\backslash\{\mu\}=\{\lambda_1,\lambda_2,\lambda_3\}=\{1,3,4\}$. The interconnection matrix $C$ represents the network topology shown in Figure \ref{fig:Ch1_ex_noiseless_topo}.
\vspace{-0.5cm}

Now, in order to design the synchronizing feedback gain for the network $\eqref{eq:Ch1_overall_dy}$ with such $C$ and the true system without using the system model, we collect input-state data as follows.
\begin{equation}
    \begin{array}{c|ccccc}
      k     \  & 0 \  & 1  \ & 2 \  & 3 \ &\ 4   \\ \hline
      x_1(k)\  & 1 \  & 1  \ & 3 \  & 4 \ &\ 5 \\
      x_2(k)\  & 0 \  &	2 \  & 1 \  & \	1  &\ 2 \\
      x_3(k)\  & 1 \  &	0 \  & 1 \  & \	2  &\ 3 \\ \hline
      u_1(k)\  & 1 \  &	-1 \ & -1 \ & \	-1 &\ {-} \\
      u_2(k)\  & -1 \ &	1 \  & 1 \  & \	1  &\ {-} \\
    \end{array}
\end{equation}
For the collected data, we observe that $ X_-$ has full row rank, but $\begin{bmatrix} X_-^\top & U_-^\top \end{bmatrix}$ does not have full row rank. This implies that $(A_\textup{true},B_\textup{true})$ cannot be uniquely identified, see \cite[Prop. 6]{van2020data}. Next, we provide the following right inverses of $X_-$:
\begin{equation}
\begin{split}
        (X^R_-)_1&=\begin{bmatrix}
           \ \frac{5}{12} \ & \frac{5}{6} \ & \frac{29}{36} \vspace{0.07cm} \ \\ 
           \ -\frac{1}{3} \ & \frac{2}{3} \ & \frac{1}{3} \vspace{0.07cm} \ \\
           \ \frac{7}{4} \ & \frac{1}{6} \ & -\frac{55}{36} \vspace{0.07cm} \ \\
           \ -\frac{13}{12} \ & -\frac{1}{2} \ & \frac{31}{36} 
        \end{bmatrix},(X^R_-)_2=\begin{bmatrix}
           \ \frac{7}{12} \ & \frac{11}{6} \ & \frac{37}{12} \vspace{0.07cm} \ \\ 
           \ -\frac{1}{3} \ & \frac{2}{3} \ & \frac{1}{3} \vspace{0.07cm} \ \\
           \ \frac{23}{12} \ & \frac{7}{6} \ & \frac{3}{4} \vspace{0.07cm} \ \\
           \ -\frac{5}{4} \ & -\frac{3}{2} \ & -\frac{17}{12} 
        \end{bmatrix},\\
        (X^R_-)_3&=\begin{bmatrix}
           \ \frac{2}{3} \ & \frac{7}{3} \ & \frac{38}{9} \vspace{0.07cm} \ \\ 
           \ -\frac{1}{3} \ & \frac{2}{3} \ & \frac{1}{3} \vspace{0.07cm} \ \\
           \ 2 \ & \frac{5}{3} \ & \frac{17}{9} \vspace{0.07cm} \ \\
           \ -\frac{4}{3} \ & -2 \ & -\frac{23}{9} 
        \end{bmatrix}.
\end{split}
\end{equation}
We see that these inverses satisfy the conditions stated in Theorem \ref{thm:Ch1_DDsyn_noiseless} (b). This, according to Theorem \ref{thm:Ch1_DDsyn_noiseless}, implies that the data $\mathcal{D}$ are informative for synchronizability with respect to $C$. Moreover, the network system \eqref{eq:Ch1_overall_dy} with $(A,B)=(A_\text{true},B_\text{true})$ is synchronized by the choice
\begin{equation}
\label{eq:Ch1_ex_gain}
    K=\frac{U_-(X^R_-)_1}{\lambda_1}=\begin{bmatrix}
        \ \frac{1}{12} \ & \frac{1}{2} \ & \frac{41}{36} \vspace{0.07cm} \ \\ 
           \ -\frac{1}{12} \ & -\frac{1}{2} \ & -\frac{41}{36} \
    \end{bmatrix}.
\end{equation}
One can verify that $A_\text{true}+\lambda B_\text{true}K$ is Schur for all $\lambda \in \sigma(C)\backslash\{\mu\}$ with such $K$.

\begin{figure}[h!]
\centering
\begin{subfigure}{0.49\textwidth}
    \includegraphics[width=0.9\linewidth]{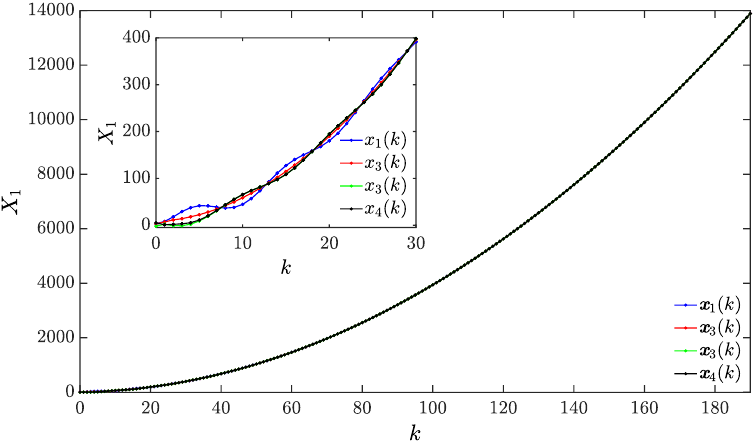}
    \caption{The trajectories on $X_1$.}
\end{subfigure}
\begin{subfigure}{0.49\textwidth}
    \includegraphics[width=0.9\linewidth]{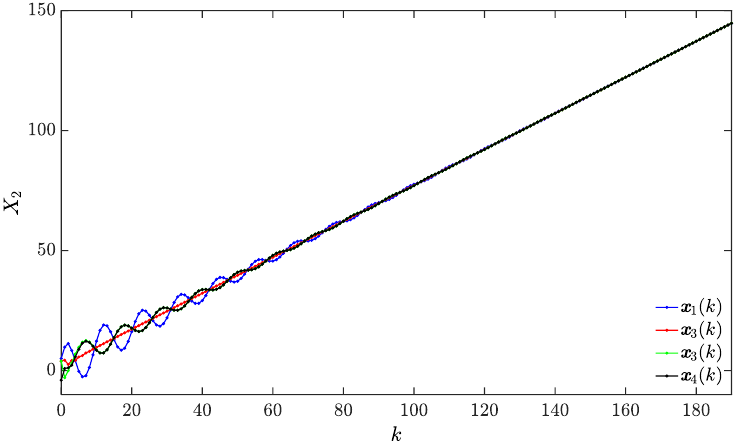}
    \caption{The trajectories on $X_2$.}
\end{subfigure}
\begin{subfigure}{0.49\textwidth}
    \includegraphics[width=0.9\linewidth]{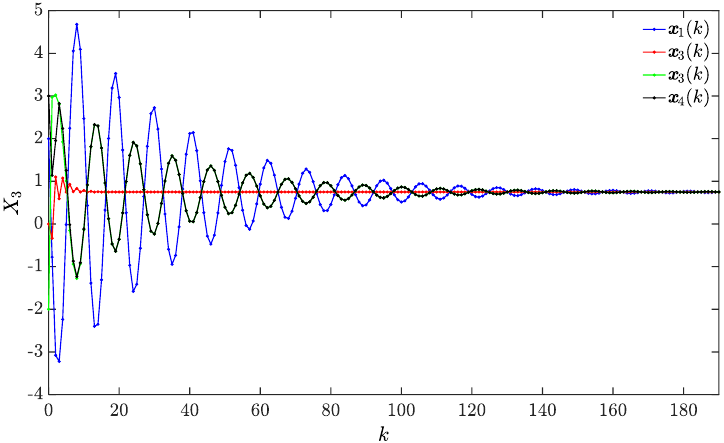}
    \caption{The trajectories on $X_3$.}
\end{subfigure}
\caption{The trajectories of systems in the network \eqref{eq:Ch1_overall_dy} on different entrances with the feedback gain in \eqref{eq:Ch1_ex_gain}.}
    \label{fig:Ch1_ex_noiseless_traj}
\end{figure}

\begin{figure}[h!]
    \centering
    \includegraphics[width=0.9\linewidth]{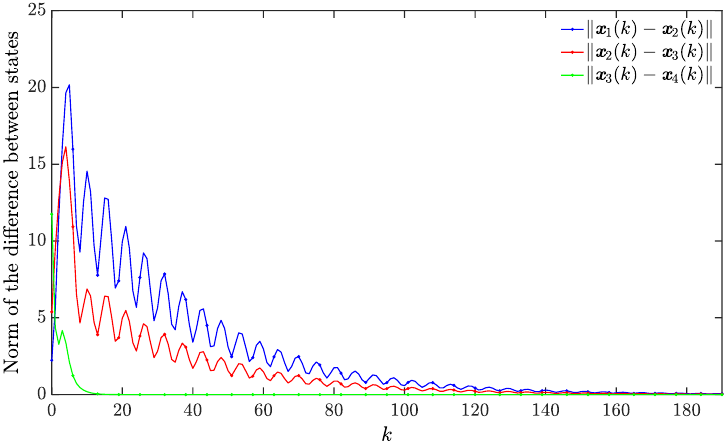}
    \caption{Norms of state differences between systems in the network \eqref{eq:Ch1_overall_dy} with the feedback gain in \eqref{eq:Ch1_ex_gain}.}
    \label{fig:Ch1_ex_noiseless_diff}
\end{figure}

Now, let the initial state of the network be
\begin{equation}
    \bm{x}=\begin{bmatrix}
        \ 3 \ & 5 \ & 2 \ & 3 \ & 4 \ & 0 \ & -2 \ & 4 \ & -2 \ & 5 \ & -4 \ & 3 \
    \end{bmatrix}^\top.
\end{equation}
Let $X_1$, $X_2$ and $X_3$ be the entrances for the states of the interconnected systems. Then, Figure \ref{fig:Ch1_ex_noiseless_traj} shows the trajectories of these systems from $k=0$ to $k=190$ on different entrances. Meanwhile, Figure \ref{fig:Ch1_ex_noiseless_diff} shows the norms of the state differences between various systems. One can observe from these figures that the systems are not stabilized but are indeed synchronized to the same state.

\section{Conclusion}

\label{Ch1_Sec_conclusion}

In this paper, we studied the data-driven synchronization problem for homogeneous network systems. For measured input-state data, we have defined an informativity notion for synchronizability with respect to a given network topology. Also, we have provided necessary and sufficient conditions for measured data to be informative for synchronizability. Finally, we have provided the corresponding design method for a feedback gain to synchronize the network directly from data. 

As a possible future work, we are interested in extending the results to a noisy data setting and unknown interconnection topologies. This means a full design of interconnection protocols directly from measured data.



\bibliography{ifacconf}             
                                                   







\end{document}